\documentclass[12pt]{amsart}
\usepackage{amscd,times,fullpage}
\usepackage[utf8x]{inputenc}

\usepackage{amsthm,amsmath}
\usepackage{amssymb}
\usepackage{graphicx}
\usepackage{amsthm}
\usepackage{enumerate}
\usepackage{dsfont}
\usepackage{mathtools}
\usepackage{tikz-cd}
\usetikzlibrary{cd}
\usetikzlibrary{positioning}
\usetikzlibrary{matrix}
\usepackage[hidelinks]{hyperref}

\hypersetup {
	pdfpagemode = {UseNone},
	pdftitle = {Ricci-iterations of well-behaved K\"ahler metrics},
	pdfauthor = {Giovanni Placini},
	pdflang = {en-UK}
} 

\newcommand{\R}{\mathbb{R}}

\newcommand{\Z}{\mathbb{Z}}
\newcommand{\T}{\mathbb{T}}

\newcommand{\CP}{\mathbb{C}\mathrm{P}}
\newcommand{\N}{\mathbb{N}}

\newcommand{\CH}{\mathbb{C}\mathrm{H}}

\newcommand{\codim}{\mathrm{codim}}

\newcommand{\Ric}{\mathrm{Ric}}

\newcommand{\C}{\mathbb{C}}            
\newcommand{\de}{\partial}

\newcommand{\K}{K\"{a}hler}

\newcommand{\ov}[1]{\overline{#1}}

\newcommand{\deb}{\ov\partial}

\newcommand{\di}{{\operatorname{d}}}

\newcommand{\vol}{{\operatorname{vol}}}

\newcommand{\Id}{\operatorname{Id}}

\newcommand{\lra}{\longrightarrow}

\newcommand{\CO}{{\mathcal{O}}}
\newcommand{\OO}{{\mathcal{O}}}

\newtheorem{theor}{Theorem}
\newtheorem{prop}{Proposition}
\newtheorem{defin}{Definition}
\newtheorem{lem}[theor]{Lemma}
\newtheorem{cor}{Corollary}

\newtheorem{ex}{Example}
\newtheorem{remark}{Remark}

\makeatletter
\newtheorem*{rep@theorem}{\rep@title}
\newcommand{\newreptheor}[2]{%
	\newenvironment{rep#1}[1]{%
		\def\rep@title{#2 \ref{##1}}%
		\begin{rep@theorem}}%
		{\end{rep@theorem}}}
\makeatother

\newreptheor{theor}{Theorem}

\usepackage{nicefrac}

\begin{document}

	\title[Ricci-iterations of well-behaved K\"ahler metrics]{Ricci-iterations of well-behaved K\"ahler metrics}
	
		\author{A.~Loi}
	\address{Dipartimento di Matematica e Informatica, Universit\'a degli studi di Cagliari, Via Ospedale 72, 09124 Cagliari, Italy}
	\email{loi@unica.it}
	\author{G.~Placini}
	\address{Dipartimento di Matematica e Informatica, Universit\'a degli studi di Cagliari, Via Ospedale 72, 09124 Cagliari, Italy}
	\email{giovanni.placini@unica.it}

	\date{\today ; {\copyright  A.~Loi and G.~Placini 2023}}
	
	\thanks{The authors are supported by INdAM and  GNSAGA - Gruppo Nazionale per le Strutture Algebriche, Geometriche e le loro Applicazioni and by GOACT - Funded by Fondazione di Sardegna.
	We acknowledge financial support by PNRR e.INS Ecosystem of Innovation for Next Generation Sardinia (CUP F53C22000430001, codice MUR ECS00000038)}
	
	\begin{abstract}
		We introduce a large class of canonical K\"ahler metrics, called in this paper well-behaved, extending metrics induced by complex space forms.
		We study K\"ahler--Ricci iterations of well-behaved metrics on compact and non-compact K\"ahler manifolds. That is, we are interested in well-behaved metrics for which the iteration of the Ricci operator is a multiple of a K\"ahler metric, i.e.,  $\rho_\omega^k=\lambda\Omega$. In particular, when $k=1$, under some condition on the maximal domain of definition of canonical coordinates, we show that $\lambda$ is forced to be positive.
		Moreover, for arbitrary $k$, we prove two additional results.
		Namely, if $\omega$ and $\Omega$ are induced by a flat metric, then $\omega$ is Ricci-flat.
		Finally, if a K\"ahler-Ricci soliton $\Omega$ arises as K\"ahler--Ricci iteration of a metric $\omega$ induced by a complex space form, then the K\"ahler--Ricci soliton is forced to be trivial, that is, K\"ahler--Einstein.
		These three theorems extend well known results on K\"ahler--Einstein metrics to higher iterations of the Ricci operator and a larger class of metrics.
\end{abstract} 
	
\maketitle
	\tableofcontents

\section{Introduction and  main results}\label{sectionint}
Yau's solution \cite{Yau78CalabiConjecture} of Calabi's conjecture is a cornerstone in the study of complex Monge-Amp\`ere equations.
It shows that, given a \K\ class $\alpha$ on a compact complex manifold $M$, any form representing $2\pi c_1(M)$ is the Ricci form of a unique \K\ form $\omega\in\alpha$. Namely, if $\rho_\omega$ denotes the Ricci form of $\omega$, the complex Monge-Amp\`ere equation
\begin{equation}\label{neweq}
	\rho_\omega=\lambda\Omega 
\end{equation}
where $[\lambda\Omega]=2\pi c_1(M)$ admits a unique solution. 
Solving similar complex Monge-Amp\`ere equations is a central problem in \K\ geometry, see \cite{Aubin78KE,ChenDonaldsonSunKEFanoI,ChenDonaldsonSunKEFanoII,ChenDonaldsonSunKEFanoIII,Tian15KStabilitySufficient} among many.
A notable case is when $\lambda\omega\in 2\pi c_1(M)$. Then equation~\eqref{neweq} is the Einstein equation
$	\rho_\omega=\lambda\omega $
and its solutions are \K--Einstein metrics. 
In this case, there are obstructions to the existence of \K--Einstein metrics when $c_1(M)>0$, cf. \cite{Futaki83Invariant,Tian97KStabilityNecessary}.
It is very natural to study solutions of equation~\eqref{neweq} and obstructions to their existence when $\omega$ and $\Omega$ satisfy additional hypotheses. As notable examples, if $\omega$ is a \K--Einstein metric induced by a non-elliptic complex space form $S$, then $M$ is a totally geodesic submanifold in $S$ \cite{Umehara87Einstein} while $\lambda>0$ if $\omega$ is a projectively induced \K--Einstein metric \cite{Hulin00KE}. 

The first part of this paper is concerned with the study of equation~\eqref{neweq} when $M$ is not necessarily compact and $\omega$ belongs to a large class of \K\ metrics including metrics induced by complex space forms or, more broadly, by generalized flag manifolds or toric manifolds.
In particular we find that, under certain conditions on the \K\ potentials of the metrics $\omega$ and $\Omega$, one necessarily has $\lambda>0$.
In order to state such hypotheses and our first result, let us recall some known facts on \K\ metrics.

Given a real analytic \K\ metric $g$ on a complex manifold $M$, one can introduce a very special \K\ potential $D^g_p$ for the metric $g$ in  a neighbourhood of a point  $p\in M$, which Calabi christened diastasis \cite{Calabi53Isometric}.
Moreover, one can always find local Bochner coordinates 
centred at $p$
such that the series expansion of the diastasis function $D^g_p$ has no terms of degree $<2$. See Section~\ref{SecWellBehaved} and references therein for details on the diastasis function and Bochner coordinates.
We can now introduce the class of \K\ metrics we are interested in, see Section~\ref{SecWellBehaved} for a detailed discussion.
\begin{defin}\label{defWB}\rm
	We call a real analytic  \K\ metric $g$ on a $n$-dimensional complex manifold $M$ {\em well-behaved at $p\in M$} if there exists  a measure zero (possibly singular or empty) real submanifold $H_p^g$ of $M$ such that the following conditions are satisfied.
	\begin{itemize}
		\item [$(A)$]
		Calabi's diastasis function $D^g_p$ around $p$ is globally defined and non-negative on $M\setminus H_p^g$;
		\item [$(B)$]
		the Bochner coordinates  $(z)=\{z_1, \dots z_n\}$ around $p$ extend to holomorphic functions $f_1, \dots , f_n$ on $M\setminus H_p^g$.
	\end{itemize}
\end{defin}

\begin{remark}\rm\label{rmkPointDependence}
	It is not clear to us whether the definition of well-behaved \K\  metric does depend on the point $p$. 
	In the compact case we even lack examples of metrics which are not well-behaved at some point, see Remark~\ref{wbforever}.
	In the non-compact case, where examples of such metrics are available (cf. Example~\ref{exnotwbnoncomp}), it is non-trivial to check whether they are well-behaved at some point.
	On the other hand, we see no reason, other than lack of examples, that suggests independence on the point.
\end{remark}

	Notice that conditions $(A)$ and $(B)$ are  always satisfied on a small enough neighbourhood $U$ of any point $p$ where Bochner coordinates are defined and $D_p^g$ is positive (except at $p$ where it vanishes) by setting  $H_p^g=\emptyset$. Thus, roughly speaking, a \K\ metric $g$ is well-behaved  at a   point  $p\in M$, if  $U$ can be enlarged to the dense open  subset  $M\setminus H_p^g$  
	in such a way that   $D_p^g$ remains non-negative and  Bochner coordinates extend to global  holomorphic functions (not necessarily coordinates).
	Obviously $H_p^g\neq\emptyset$ if $M$ is compact.

In contrast with this, we will require a global condition on the metrics. Such condition was first considered in \cite{ArezzoLoi04EuclideanVolume}, see also \cite{ArezzoLoiZuddas16Toric}. 
\begin{defin}\label{defvoleucl}\rm
Let $M$ be a $n$-dimensional complex manifold and  $g$ be a well-behaved \K\ metric at a point $p\in M$. 
We define the {\em Bochner-Euclidean volume} of $g$ at $p$ to be
\begin{equation}\label{voleucl}
\vol_{eucl} (M\setminus H_p^g, g):=\int_{M\setminus H_p^g}d\mu(f)
\end{equation}
where $f_1, \dots , f_n$ are as in Definition~\ref{defWB} and $d\mu (f):=\frac{i^n}{2^n}df_1\wedge d\bar f_1\wedge\cdots df_n\wedge d\bar f_n$.
\end{defin}
Notice that  $\vol_{eucl} (M\setminus H_p^g, g)$
strongly depends on the Bochner coordinates
and so on the metric $g$.
One can think of the Bochner--Euclidean volume 
as the volume of $M\setminus H_p^g$
with respect of  the pull-back of the standard Lebesgue measure of $\C^n$
via  the map $M\setminus H_p^g\rightarrow \C^n, q\mapsto (f_1(q), \dots , f_n(q))$.

We are now ready to state our first result, Theorem~\ref{mainteor}. 
Although its hypotheses appear technical at first glance, they are satisfied by a large class of \K\ metrics, see Section~\ref{SecWellBehaved} for several examples.
In fact, when specialized to certain canonical metrics, Theorem~\ref{mainteor} yields several corollaries (Corollary~\ref{Cor1}-\ref{Cor3} below) which extend well known results.
\begin{theor}\label{mainteor}
Let $M$ be a complex $n$-dimensional manifold and let  $g, G$ be two well-behaved   \K\ metrics at the same point  $p\in M$ satisfying
\begin{equation*}
\rho_\omega=\lambda\Omega,
\end{equation*}
where $\omega$ and $\Omega$ are the \K\ form associated 
to $g$ and $G$ respectively.
Assume that $g$ and $G$ satisfy the following conditions:
\begin{enumerate}
\item [(i)]
$\vol (M, \omega)=\int_M\frac{\omega^n}{n!}<\infty$;
\item [(ii)]
$\vol_{eucl} (M\setminus H_p^g, g)=\infty$;
\item [(iii)]
$\codim_{\R}(H_p^g)\geq 2$.
\end{enumerate}
with $H_p^g$  as in Definition~\ref{defWB} .
Assume that also  $G$ is well-behaved at $p$ and 
\begin{enumerate}
\item [(iv)]
$\codim_{\R}(H_p^G)\geq 2.$
\end{enumerate}

Then  
$\lambda >0$.
\end{theor}

The simplest case, which motivated Theorem~\ref{mainteor}, is that of compact homogeneous \K\ manifolds.
Indeed, the Ricci tensor of a homogeneous \K\ metric is always a homogeneous \K--Einstein metric, cf. Proposition~\ref{PropHomogeneous} below. 
If in addition the homogeneous \K\ form is integral, it is projectively induced (see \cite{Takeuchi78Homogeneous,LoiZuddasGromovWidthHKM}).
To the best of our knowledge, this is the only example of projectively induced \K\ forms $\omega$, $\Omega$ satisfying the hypotheses of Theorem~\ref{mainteor}.
Inspired by this case we study projectively induced (radial) metrics $g$ on $\CP^1$ and show in Proposition~\ref{PropProjIndCP1} that their Ricci tensor is projectively induced if and only if it is \K--Einstein.
Moreover, we give a very explicit example: Ricci forms of well-behaved metrics on $\CP^1$, see Example~\ref{ExPIntoWB}. 
In particular, the positivity of the first Chern class allows us to construct well-behaved, not projectively induced metrics $g$ and $G$ on $\CP^1$ which satisfy the hypotheses of Theorem~\ref{mainteor}.

Some comments on the necessity of the hypotheses (i)-(iii) are in order.

\noindent - Assumption (i) is a necessary condition for Theorem~\ref{mainteor}
to be true. For example the hyperbolic metric $g_{hyp}$ on the complex hyperbolic space is well-behaved  and satisfies conditions (ii) and (iii) (see Example~\ref{hypmet}) and being \K-Einstein  
it satisfies \eqref{neweq} with $g=G$ and 
$\lambda=-2(n+1)<0$.

\noindent - Also assumption (ii) is necessary. Take for example any bounded domain in $\C^n$
with the flat metric. Examples of compact manifolds not satisfying (ii) nor (iii) are complex tori with the flat metric (Example~\ref{exnotwbcomp}) or Riemann surfaces with the hyperbolic metric and products thereof (Example~\ref{exnotwbcomp2}).

\noindent - Concerning the necessity of property (iii), we do not know of any compact \K\ manifold satisfying (i) and (ii) but not (iii), see Examples~\ref{exnotwbcomp} and~\ref{exnotwbcomp2} below. In particular, it would be interesting to understand the interplay between conditions (ii) and (iii).

Notice that the assumptions (i)-(iii) of Theorem~\ref{mainteor} are satisfied, for instance, by the complex projective space $\CP^n$.
Moreover, properties (i)-(iii) are hereditary for projective varieties or, more generally, compact \K\ submanifolds of certain compactifications of $\C^n$, see Proposition~\ref{proptoric}.
From this we obtain the following results.

\begin{cor}\label{Cor1}
Let $M$ be a compact complex manifold endowed with two projectively induced \K\ forms $\omega$, $\Omega$ satisfying
$\rho_\omega=\lambda\Omega.$ Then $M$ is Fano.
\end{cor}
\begin{cor}\label{Cor2}
	Let $M$ be a compact complex manifold equipped with two \K\ forms $\omega$, $\Omega$ induced by immersions into classical flag manifolds. If $\rho_\omega=\lambda\Omega$, then $M$ is Fano.
\end{cor}
Observe that Corollary~\ref{Cor1} is in fact implied by Corollary~\ref{Cor2}. On the other hand, not all classical flag manifolds are projectively induced, even among irreducible ones and up to homothety, cf. \cite{LoiMossa23Rigidity}.

\begin{cor}\label{Cor3}
	Let $M$ be a compact complex manifold equipped with two \K\ forms $\omega$, $\Omega$ induced by immersions $\varphi_1,\ \varphi_2$ into a toric manifold such that $\varphi_1(p)=\varphi_2(p)$ is fixed by the torus action for some $p\in M$. If $\rho_\omega=\lambda\Omega$, then $M$ is Fano.
\end{cor}

As mentioned above, the Einstein equation $\rho_\omega=\lambda\omega$ is a special case of  equation~\eqref{neweq}.
In light of this, if $\omega=\Omega$ in  Corollary~\ref{Cor1} and Corollary~\ref{Cor3} we recover some known results on \K--Einstein submanifolds.
Namely, the main theorem of \cite{Hulin00KE} and \cite[Theorem~2.6]{ArezzoLoiZuddas16Toric} respectively. 
In fact, the proof of Theorem~\ref{mainteor} is inspired by these results and their proofs.

One also derives the following consequence for a  \K\--Einstein metric $g$  on a compact  complex manifold $M$. 
If $c_1(M)<0$ or $c_1(M)=0$, the existence of \K\--Einstein metrics is guaranteed by \cite{Aubin78KE} and \cite{Yau78CalabiConjecture}, but explicit examples are extremely rare. 
By Theorem~\ref{mainteor} such metrics necessarily satisfy at least one of the following conditions for a point $p\in M$: 
\begin{itemize}
	\item 
	$g$ is  not well-behaved at $p$;
	\item
	$\vol_{eucl} (M\setminus H_p^g, g)<\infty$;
	\item
	$\codim_{\R}(H_p^g)< 2$.
\end{itemize}

Let us  now comment on the necessity of the condition of well-behaviour of the metrics $g$ and $G$ in Theorem~\ref{mainteor}.
We observe that requiring both  $g$ and $G$ to be well-behaved is necessary in the compact case.
For instance, one can consider a compact complex manifold with $c_1(M)<0$ and a \K\ form $\Omega$ such that $\lambda\Omega$ represents the first Chern class (up to a $\frac{i}{2\pi}$ factor).
By Calabi's conjecture, in a given \K\ class there exists a unique \K\ form $\omega$ whose Ricci form $\rho_\omega$ equals $\lambda\Omega$.
By Tian's theorem \cite{Tian90Approximation}, if we choose the class $[\omega]$ to be integral, we can approximate the form $\omega$ with (multiples of) projectively induced metrics $\omega_k$.
For $k$ large enough $\rho_{\omega_k}/\lambda$ is still positive, hence a \K\ form. 
We then have projectively induced \K\ forms $\omega_k$ satisfying
$$\rho_{\omega_k}=\lambda\Omega_k$$
for \K\ forms $\Omega_k$ which cannot be well-behaved by Theorem~\ref{mainteor}.
In fact one can see that the necessary condition is weaker than $G$ being well-behaved, cf. Corollary~\ref{CorConditionsProjInduced}.

Analogously, if we pick $\lambda\in2\pi\Z$ we can find projectively induced metrics $\Omega_k$ satisfying
$$\rho_{\omega_k}=\lambda\Omega_k$$
for \K\ forms $\omega_k$ which cannot be well-behaved by Theorem~\ref{mainteor}.

	\vspace{0.3cm}
	
In the second part of this paper we focus on a generalization of equation~\eqref{neweq} related to the so-called \textit{\K--Ricci iterations}. Namely, if $\pm \rho_\omega$ is again a  \K\ metric we can consider its Ricci form $\rho^2_\omega$ or, more generally, we can define inductively
\begin{equation}\label{EqDefIteration}
\rho^k_\omega:=\begin{cases}
	\omega, & \text{if } k=0\\
	\rho_{\alpha_{k-1}}, & \text{if $\alpha_{l}:=\pm\rho^{l}_\omega>0$, for all $l<k$}.
\end{cases}
\end{equation}
This iteration was introduced by Nadel in \cite{Nadel95RicciIterations} where he also proved that periodic points of order two or three must be K\"ahler-Einstein metrics.
This was generalized \cite{Keller09RicciIterations,Rubinstein08ConvergenceRicciIteration} to periodic points $\omega$ of any order, that is, those satisfying $\rho^k_\omega=\lambda\omega$ for any $k\in\N$.
The iteration \eqref{EqDefIteration} can be regarded as a discretization of the \K--Ricci flow.
Observe that, if the manifold $M$ is compact, by Calabi's conjecture one does not need positivity assumptions to reverse the construction above and define the  \textit{inverse Ricci--\K\ iterations} $\rho^{-k}_\omega$.
Both of this discrete dynamical systems have been studied in the literature, see for instance \cite{Berman19ConvergenceRicciIteration,DarvasRubinstein19ConvergenceRicciIterations}.

Here we focus on the following generalized Monge-Amp\`ere equation on a complex manifold $M$
\begin{equation}\label{EqRicciIteration}
	\rho^k_\omega=\lambda\Omega
\end{equation}
 where $\Omega$ is again a \K\ form.
In particular we study equation~\eqref{EqRicciIteration} for a special class of well-behaved \K\ metrics.
Our second result is the following theorem dealing with \K\ metrics induced by the flat metric (which are well-behaved by Example~\ref{flatmetric} below).

\begin{theor}\label{mainteor2}
Let $M$ be a  complex  manifold with two metrics $g$ and $G$ induced by the flat metric.
If the corresponding \K\ forms $\omega,\Omega$ satisfy 
$\rho_\omega^k=\lambda\Omega$ for some $k\geq 1$ and $\lambda\in\R$,
then $(M,g)$ is a totally geodesic submanifold of the flat ambient space.
\end{theor}

Observe that this result can be seen as a generalization of \cite[Theorem~2.1]{Umehara87Einstein}. 
Namely, when $k=1$ and $g=G$, Theorem~\ref{mainteor2} shows that a \K--Einstein submanifolds of $\C^n$ with the flat metric is necessarily Ricci-flat, hence a totally geodesic submanifold.
It is worth pointing out that the metric $g$ is assumed to be induced by a flat one for simplicity, but the hypothesis on $g$ can be sensibly relaxed, cf. Remark~\ref{RmkAlsoOtherMetrics} below.

Our third and last result deals with \K--Ricci solitons (KRS). Recall that a \K\ metric $G$ is a KRS  if there exists a holomorphic vector field $X$ (the solitonic vector field) such that $\Ric(G)=\mu G+L_XG$ where $\Ric(G)$ denotes the Ricci tensor of $G$ and $L_X$ the Lie derivative in the direction of $X$.
\K--Ricci solitons are important generalizations of \K--Einstein metrics which arise in the study of the \K--Ricci flow. 
Our next result show that KRS cannot arise as \K--Ricci iterations of metrics induced by complex space forms.

\begin{theor}\label{mainteor3}
Let $M$ be a  complex  manifold with two real analytic \K\ metrics $g$ and $G$.
Assume that $g$ is induced by a complex space form and $G$ is a KRS.
If the corresponding \K\ forms $\omega,\Omega$ satisfy 
$\rho_\omega^k=\lambda\Omega$ for some $\lambda\neq 0$ and $k\geq 0$,
then $G$ is trivial, i.e. \K--Einstein.
\end{theor}
Observe that, when $\lambda=0$, one cannot draw any conclusion. Take for instance $g$ to be the flat metric on $\C$ and $G$ to be Hamilton's cigar KRS \cite{Hamilton88RicciFlow}.
It is worth mentioning that we do not know of any \K--Einstein metric $G$ and \K\ metric $g$ induced by a complex space form which satisfy $\rho_\omega^k=\lambda\Omega$ for $\lambda<0$ unless $g=G$ and $k=1$.
 In that case $M$ is forced to be a totally geodesic submanifold in $\CH^n$ \cite{Umehara87Einstein}.
 On the other hand, such examples for $\lambda>0$ are discussed in Section~\ref{SecProofs}, see in particular Proposition~\ref{PropHomogeneous}.

Also in this case we generalize some recent results on KRS induced by complex space forms which indeed served as partial motivation for this paper. 
In particular, if we assume $k=0$ (and $\lambda>0$) in Theorem~\ref{mainteor3}, we recover \cite[Theorem~1.1]{LoiMossa21KRSIntoComplexSpaceForms}. 
If we restrict to compact complex manifolds, then this result can be rephrased in terms of the inverse \K--Ricci iterations. 
Namely, Theorem~\ref{mainteor3} shows that none of the inverse \K--Ricci iterations $\rho^{-k}_\omega$ of a non-trivial KRS $(g,X)$ on a compact complex manifold $M$ can be induced by a complex space form.

Finally, observe that the metric $g$ is assumed to be induced by a complex space form, but the proof of the theorem shows that this hypothesis can be relaxed, cf. Remark~\ref{RmkAlsoOtherMetrics} below.
\vskip 0.2 cm

\noindent\textbf{Organization of the paper.} In Section~\ref{SecWellBehaved} we discuss well-behaved \K\ metrics as well as the conditions (i)-(iii) of Theorem~\ref{mainteor}. 
Moreover, we prove that metrics induced by suitable compactifications of $\C^n$ are both well-behaved and satisfy conditions (i)-(iii), see Proposition~\ref{proptoric}. 
This is crucial in the proofs of Corollaries~\ref{Cor1}-\ref{Cor3} which complete Section~\ref{SecWellBehaved}.
In Section~\ref{SecProofs} we give the proof of Theorem\ref{mainteor} and a discussion on \K--Ricci iterations on compact homogenous \K\ manifolds. 
In particular, we give a very explicit example:  \K--Ricci iterations of well-behaved metrics on $\CP^1$. 
Finally, Section~\ref{SecFinal} is devoted to the proofs of Theorem~\ref{mainteor2} and Theorem~\ref{mainteor3}.

\section{Well-behaved \K\ metrics}\label{SecWellBehaved}
In order to discuss well-behaved metrics, we  briefly recall  some basic facts on Calabi's diastasis functions and Bochner coordinates.
Given a complex manifold $M$ endowed with a real analytic \K\ metric $g$, one can introduce on  a neighbourhood of a point $p\in M$, a very special \K\ potential $D^g_p$ for the metric $g$, which Calabi christened {\em diastasis}.
Recall that a \K\ potential is an analytic function $\Phi$ defined in a neighbourhood of a point $p$ such that $\omega =\frac{i}{2}\partial \bar\partial\Phi$, where $\omega$ is the \K\ form associated to $g$.
In a complex coordinate system $(z)=\{z_1, \dots ,z_n\}$ around $p$ one has
\begin{equation}\label{locg}
	g_{\alpha\bar\beta}=
	2g(\frac{\partial}{\partial z_{\alpha}},
	\frac{\partial}{\partial \bar z_{\beta}})
	=\frac{{\partial}^2\Phi}
	{\partial z_{\alpha}\partial\bar z_{\beta}}.
\end{equation}
A \K\ potential is not unique: it is defined up to summing with the real part of a holomorphic function.
By duplicating the variables $z$ and $\bar z$, a potential $\Phi$ can be complex analytically
continued to a function  $\tilde\Phi$ defined on a neighbourhood $U\subset M\times\overline M$ of the diagonal containing $(p, \bar p)$ (here $\overline M$ denotes the conjugated manifold of $M$).
The {\em diastasis function} is the \K\ potential $D^g_p$ around $p$ defined by
$$D^g_p(q)=\tilde\Phi (q, \bar q)+
\tilde\Phi (p, \bar p)-\tilde\Phi (p, \bar q)-
\tilde\Phi (q, \bar p).$$
The diastasis function is characterized among potentials by the property that in every coordinates system  $(z)$ centred in $p$
$$D^g_p(z, \bar z)=\sum _{|j|, |k|\geq 0}
a_{jk}z^j\bar z^k,$$
with  $a_{j 0}=a_{0 j}=0$ for all multi-indices $j$.
More generally, for later use, we give the following (see also \cite{LoiMossa21KRSIntoComplexSpaceForms,LoiMossa22HBD})
\begin{defin}\label{defdiast}\rm
	Let $M$ be a complex manifold and $f:U\rightarrow \R$ a real analytic function defined on a neighbourhood   $U$ of a point $p\in M$.
	We say  that $f$  is of  {\em diastasis type} if 
	in one (and hence any) coordinate system 
	$\{z_1, \dots , z_n\}$ centred at $p$ the   expansion of $f$  in $z$ and $\bar z$ 
	does not contain non-constant purely holomorphic or anti-holomorphic terms (i.e. of the form $z^{j}$ or $\bar{z}^{j}$ with  $\vert j\vert > 0$). 
\end{defin}

\begin{remark}\label{remdefdiast}\rm
	By its very definition Calabi's diastasis function $D_p^g$ is of diastasis type.
\end{remark}

\begin{remark}\label{locexp(1,1)}\rm
	Calabi's procedure can 
	be applied to any (not necessarily \K) real-analytic real-valued  closed form of type $(1, 1)$. Namely, if $\Gamma$ is a such a form, on a neighbourhood of a point $p\in M$ one can write
	$\Gamma =\frac{i}{2}\partial\bar\partial \gamma$,  with $\gamma$ of diastasis type.
\end{remark}
One can always find  local complex coordinates  $(z)$ centred at $p$ such that
$$D^g_p(z, \bar z)=|z|^2+\sum _{|j|, |k|\geq 2} b_{jk}z^j\bar z^k.$$
These coordinates, called {\em Bochner coordinates}
around $p$, are 
uniquely defined up
to a unitary transformation
(cf. \cite{Bochner47coordinates,Calabi53Isometric} or \cite{Loi18Book} for an updated account on the subject).
The following result  summarizes some fundamental properties of Calabi's diastasis function and Bochner coordinates.

 \begin{theor}[\cite{Calabi53Isometric}]\label{ThmCalabi}
  Let $\varphi :(S, h)\rightarrow (M, g)$ be a holomorphic isometric immersion between \K\ manifolds of complex dimension $m$ and $n$ respectively.
 If $g$ is real analytic, then $h$ is real analytic and for every point $q\in S$ 
 $$D^h_q=\varphi ^*D^g_{\varphi (q)}.$$ 
 Moreover, if $(w_1,\dots ,w_m)$ is a system of Bochner coordinates in a neighbourhood $U$ of $q$, then there exists a system of Bochner coordinates $(z_1,\dots ,z_n)$ centred at $\varphi (q)$ such that
 \begin{equation}\label{zandZ}
 	z_1\circ \varphi |_{U}=w_1,\dots , z_m\circ \varphi |_{U}=w_m.
 \end{equation}
 \end{theor}
We provide now several examples of well-behaved \K\ metrics both on non-compact and compact complex manifolds aimed at showing that they are fairly common among real analytic \K\ metrics.

\begin{ex}\rm (Products and linear combinations)\label{opwb}
	Let $M$ be any complex manifold.
	It is not hard to verify that if $g_j$, $j=1, \dots ,k$, are well-behaved \K\ metrics at a point $p\in M$ and  $\alpha_j$ are positive real  numbers such that $g=\sum_{j=1}^k\alpha_jg_j$ is still a \K\ metric,  then $g$ satisfies property (A) in Definition~\ref{defWB}. 
	On the other hand, we cannot describe the Bochner coordinates of $g$ in terms of those of the metrics $g_j$. Thus, we do not have any information on whether linear combinations of well-behaved metrics satisfy condition (B) in Definition~\ref{defWB}. 
	
	A  different behaviour is exhibited by \K\ products.
	In fact, if $(M_j, g_j)$ with $j=1, \dots ,k$ are \K\ manifolds with $g_j$
	well-behaved at $p_j\in M_j$ for all $j$ , then the \K\ metric $g_1\oplus\cdots\oplus g_k$ 
	on $M_1\times\cdots \times M_k$ is  well-behaved at $(p_1, \dots, p_k)$. 
	Moreover, it is easy to see that products preserve assumptions (i)-(iii) in Theorem~\ref{mainteor}.
\end{ex}

\begin{ex}\rm (Hereditary property)\label{hered}
	By  Theorem~\ref{ThmCalabi}  we get many examples of well-behaved \K\ metrics induced by well-behaved ones.  More precisely, suppose that $\varphi: S\rightarrow M$ is a holomorphic isometric immersion  between \K\ manifolds $(S,  h)$ and $(M, g)$. If $g$ is well-behaved at $p=\varphi (q)$ for some $q\in S$,  then $h$ is well-behaved at $q$. Moreover, we have $H_q^h=\varphi^{-1}(H^g_{p})$.
	
	It is easy to exhibit examples showing that conditions (i) and (ii) are not hereditary, see also Remark~\ref{RmkEsempioCinCPn}. 
	Therefore, if $(M, g)$ is a \K\ manifold satisfying conditions (i)-(iii) in Theorem~\ref{mainteor}, we cannot conclude that a \K\ submanifold $(S,  h)$ inherits those properties.
\end{ex}

\begin{ex}\rm\label{hypmet} (Hyperbolic disc)
	The hyperbolic metric $g_{hyp}$ on $\CH^n=\{z=(z_1, \dots z_n) \ | \ \in \C^n \ |\ |z|^2<1\}$ is well-behaved at 
	any point of $\C H^n$.
	This can be verified at the origin $0\in \C H^n$.
	Indeed, the  \K\ form  associated to $g_{hyp}$ is given by $\omega_{hyp}=-\frac{i}{2}\partial\bar\partial\log (1-|z|^2)$ whose diastasis at the origin $0\in\C^n$ reads as $D_0^{g_{hyp}}(z)=-\log (1-|z|^2)$. 
	Hence, the Bochner coordinates are the restriction of the Euclidean coordinates to $\C H^n$ and $H_0^{g_{hyp}}=\emptyset$. 
This shows that $g_{hyp}$ is well-behaved and that its Euclidean volume $\vol_{eucl}(\CH^n,g_{hyp})$ is finite although the volume $\vol(\CH^n,g_{hyp})$ is infinite.
	More generally, one can find  in \cite{DiScalaLoi07SymmetricSpaces} a description of Calabi's diastasis function and Bochner coordinates for bounded symmetric domains with their Bergman metrics which yields the same conclusion.
\end{ex}

\begin{ex}\rm\label{flatmetric} (Flat $\C^n$)
	The flat metric $g_0$ on $\C^n$ is well-behaved at any point of $\C ^n$.
	This too can be easily verified at the origin $0\in \C^n$.
	Indeed, the \K\ form  associated to $g_{0}$ is given by $\omega_{0}=\frac{i}{2}\partial\bar\partial |z|^2$
	whose diastasis at the origin reads as $D_0^{g_{0}}(z)=|z|^2$. Thus, the Bochner coordinates
	are the Euclidean coordinates  and $H_0^{g_0}=\emptyset$.
	It follows by Example~\ref{hered} that the \K\ metric $g$ induced on a Stein manifold $M$ by the flat metric
	through an  embedding into a complex Euclidean space is  indeed well-behaved at any $p\in M$ with $H_p^g=\emptyset$.
	Clearly, the Bochner-Euclidean volume coincides with the \K\ volume in this case and they are infinite.
	
\end{ex}
\begin{ex}\rm\label{projspace} (Fubini-Study metric)
	Let $g_{FS}$ be the Fubini-Study metric of $\C P^n$ whose associated \K\ form is given in homogeneous coordinates by $\omega_{FS}=\frac{i}{2}\partial\bar\partial\log (|Z_0|^2+\cdots +|Z_n|^2)$.
	It is easy to see that  the diastasis at the point $p_0=[1, 0, \dots , 0]$ is globally defined on $\C P^n\setminus H^{g_{FS}}_p$, where 
	$H_{p_0}^{g_{FS}}$ is the hyperplane $\{Z_0=0\}$.
	Moreover,  the  affine coordinates $z_j=\frac{Z_j}{Z_0}$, $j=1, \dots n$ are indeed global  Bochner coordinates on $\C P^n\setminus H^{g_{FS}}_p$ because the diastasis in these coordinates is given by
	$D_{p_0}^{g_{FS}}(z)=\log (1+|z|^2)$.
	Hence $g_{FS}$ is well-behaved at $p_0$ (in fact at all points because $\omega_{FS}$ is homogeneous).
	It is then clear that the Bochner-Euclidean volume $\vol_{eucl}(\CP^n\setminus H^{g_{FS}}_p,g_{FS})$ is infinite while $\vol(\CP^n,g_{FS})$ is clearly finite.
	Thus $(\CP^n,g_{FS})$ is the simplest example of compact \K\ manifold with a well-behaved metric which satisfies (i)-(iii) of Theorem~\ref{mainteor}.
\end{ex}

\begin{ex}\rm (Projectively induced metrics)\label{projindnoncomp}\label{projindwb}
	A \K\ metric $g$ on a (not necessarily compact) complex manifold $M$ is said to be  {\em projectively induced} if there exists a holomorphic immersion $\varphi: (M, g)\rightarrow (\C P^{N\leq \infty}, g_{FS})$ such that $\varphi^*g_{FS}=g$ (see \cite{Calabi53Isometric}). 
	Here $\C P^{N\leq \infty}$ denotes the  finite or infinite dimensional complex projective space equipped with the Fubini-Study metric $g_{FS}$. 
	It  is not hard to see that the second part of Calabi Theorem~\ref{ThmCalabi} still holds when the ambient space is  the infinite  dimensional complex projective space $(\C P^{\infty}, g_{FS})$.
	Therefore, as in Example~\ref{hered}, one sees that a projectively induced metric $g$ on a complex manifold $M$ is well-behaved at any point $p\in M$. 
	Many examples of compact and non-compact projectively induced \K\ metrics have been considered in the literature (see \cite{Loi18Book} and references therein).
	Moreover, when $N<\infty$ and $M$ is compact, a projectively induced metric $g$ on $M$ has infinite Bochner--Euclidean volume (and finite volume), cf. Proposition~\ref{proptoric} below.
\end{ex}

In the following we deal with  two examples of  
well-behaved \K\  metrics on compactifications of $\C^n$ which are not necessarily projectively induced.
\begin{ex}\rm\label{genflag}(Generalized flag manifolds)
	Let $(M, g)$ be a generalized flag manifold, that is, a compact and simply-connected homogeneous \K\  manifold.
	If  the \K\ form associated to $g$ (respectively one of  its multiples) is  integral  then the  \K\ metric $g$ (resp. its integral multiple)  is projectively induced  (see \cite{LoiZuddasGromovWidthHKM})
	and hence well-behaved at any point  by Example~\ref{projindwb}. 
	Typical examples are homogeneous \K\ metrics of Hermitian symmetric spaces of compact type or \K-Einstein metrics on generalized flag manifolds.
	On the other hand, there exist   generalized flag manifolds $(M, g)$
	where $g$ and all its multiples are not integral, but 
	$g$ is well-behaved. For example $(\C P^1\times \C P^1, g)$, where $g=g_{FS}\oplus\sqrt{2}g_{FS}$  is well-behaved  by Example~\ref{opwb}.
	
	More generally, there exists a large class of  irreducible classical generalized flag manifolds $(M, g)$ with $g$ well-behaved and not necessarily projectively induced (see \cite[Theorem~1]{LoiMossaZuddas19BochnerFlag}).
	Such manifolds are compactifications of $\C^n$, namely   $M=\C^n \sqcup H$, where $H$ is a complex (possibly singular) submanifold of complex codimension one.
	Moreover, the diastasis function $D_0^g$ at $0\in\C^n$ is defined on $\C^n$ and positive away from $0$  and the Bochner coordinates are globally defined on $\C^n$. This shows that the Euclidean volume $\vol_{eucl}(M\setminus H,g)$ is infinite while, by compactness $\vol(M,g)<\infty$. Therefore, such classical generalized flag manifolds satisfy conditions (i)-(iii) of Theorem ~\ref{mainteor}.
	
	Notice that for  the Hermitian symmetric spaces of compact type
	$H_p^g$ is in fact the cut locus with respect to the point $p$ (see for instance \cite{Loi06DiastasisHSS,LoiMossa11DiastaticExp,Tasaki85CutLocusDiastasis}).
	We believe (see also \cite[Conjecture p.~499]{LoiMossaZuddas19BochnerFlag}) that this property should characterize Hermitian symmetric spaces among generalized flag manifolds.
\end{ex}

\begin{ex}\rm\label{toric}(Toric manifolds)
	Recall that a toric manifold $M$ is a complex manifold which contains an open dense subset biholomorphic to $(\C^*)^n$ and such that the canonical action of $(\C^*)^n$ on itself 
	extends to a holomorphic action on the whole manifold $M$.
	A toric K\"ahler metric $\omega$ on $M$ is a K\"ahler metric which is invariant under the action of the real torus $T^n = \{ (e^{i \theta_1}, \dots , e^{i \theta_n}) \ | \ \theta_i \in \R \}$ contained in the complex torus $({\C}^*)^n$.
	That is, for every fixed $\theta \in T^n$ the diffeomorphism $f_{\theta}: M \rightarrow M$ given by the action of $(e^{i \theta_1}, \dots , e^{i \theta_n})$ is an isometry.
	
	We have the following well known fact on toric manifolds, see for example \cite[Section~2.2.1]{Donaldson04ConjInKahler}, \cite[Proposition~2.18]{Batyrev93QuantumCohomology} or refer to  \cite[Appendix]{ArezzoLoiZuddas16Toric}  for a self-contained proof.
		If $M$ is a projective, compact toric manifold, then there exists an open dense subset $X \subseteq M$ which is algebraically biholomorphic to $\C^n$. More precisely, for every point $p \in M$ fixed by the torus action there are an open dense neighbourhood $X_p$ of $p$ and a biholomorphism $\phi: X_p \rightarrow \C^n$ such that $p$ is sent to the origin and the restriction of the torus action to $X_p$ corresponds via $\phi$ to the canonical action of $(\C^*)^n$ on $\C^n$. Moreover, one can show \cite[Theorem~2.6]{ArezzoLoiZuddas16Toric} that the diastasis function $D_p^g$ of the toric metric $g$ is defined and positive on $X_p$ and that the coordinates on $X_p$ given by $\phi$ are Bochner coordinates for $g$.
		 Therefore invariant metrics on toric manifolds are well-behaved and have infinite Bochner--Euclidean volume. Since they have finite volume by compactness, toric metrics satisfy conditions (i)-(iii) of Theorem~\ref{mainteor}.
\end{ex}

\begin{ex}\rm\label{exnotwbcomp}(Flat  torus)
	Let $\T^n={\C}^n/{\Z}^{2n}$ be the $n$-dimensional complex torus endowed with the flat metric $g_0$.
	Let
	$$M_0=\{(z_1, \dots z_n)\in\C^n \ | \  -\frac{1}{2}<\Im (z_j)<\frac{1}{2} \ , \ -\frac{1}{2}<\Re (z_j)<\frac{1}{2} \ , \   j=1, \dots ,n\}$$  
	be the fundamental domain of the $\Z^{2n}$-action containing the origin  $0\in\C^n$ where $\Re (z)$ and $\Im (z)$ are the real and imaginary part of $z\in\C$ respectively.
	The diastasis $D_0^{g_0}$ of $g_0$ with respect to $0$ reads as $D^{g_0}_{0}(z, \bar z)=|z|^2=\sum _{j=1}^{n}|z_j|^2$ and the Euclidean coordinates $(z_1,\dots ,z_n)$ are in fact Bochner coordinates for $g_0$ around $0$.
	Notice that $M_0=M\setminus H^{g_0}_0$ is the maximal domain of extension of $D^{g_0}_0$ so that the Bochner--Euclidean volume $\vol_{eucl}(M\setminus H^{g_0}_0,g_0)$ is finite.
	Moreover, $H^{g_0}_0$ is given by 
	$$H_0^{g_0}=\bigcup_{j=1}^n\left\{(z_1, \dots z_n)\in \T^n \ |\ \Re (z_j)=\frac{1}{2}\right\}\ \bigcup\ \bigcup_{j=1}^n\left\{(z_1, \dots z_n)\in \T^n \ |\ \Im (z_j)=\frac{1}{2}\right\}$$
	which is a singular submanifold of $M$ of real codimension $1$.
	Hence  $g_0$ is well-behaved but not projectively induced  by Example~\ref{projindwb} (or even locally projectively induced by \cite{Calabi53Isometric}).
\end{ex}

\begin{ex}\rm\label{exnotwbcomp2} (Riemann surfaces)
	The hyperbolic metric $g_{hyp}$ on a compact Riemann surface  $\Sigma_g$ of genus $g\geq 2$ is well-behaved.
	Indeed, consider $\Sigma_g=\CH^1/\Gamma_g$ where $\Gamma_g=\pi_1(\Sigma_g)$ and let $D_0^{g_{hyp}}$ be the diastasis function at $0\in\CH^1$.
	Then its maximal domain of definition $M_0$ is the fundamental domain of the $\Gamma_g$-action and satisfies $M_0=\Sigma_g\setminus H_0^{g_{hyp}}$ where  $H_0^{g_{hyp}}$ is a bouquet of $2g$ circles.
	Also in this case the euclidean coordinate $z$ is the  Bochner coordinate for $g_{hyp}$ around $0$.
	Hence $g_{hyp}$ is well-behaved, but does not satisfy condition (ii) in Theorem~\ref{mainteor} because its Euclidean volume is finite. 
	Moreover, $g_{hyp}$ does not satisfy property (iii) because $H_0^{g_{hyp}}$ has real codimension $1$ in $\Sigma_g$.
	Hence  $g_{hyp}$ is well-behaved but not projectively induced by Example~\ref{hypmet} (or even locally by \cite{Calabi53Isometric}). 

	Similar arguments apply to a compact quotient 
	of $\C H^n$ equipped with the hyperbolic metric such as products of compact Riemann surface  $\Sigma_{g_j}$ of genus $g_j\geq 2$.
\end{ex}

Now  we exhibit the existence of not well-behaved \K\ metrics on non-compact manifolds.

\begin{ex}\rm(Metrics which are not well-behaved at a point)\label{exnotwbnoncomp} 
	A necessary condition for $g$ to satisfy property $(A)$ of Definition~\ref{defWB} at a point $p$ is that the   maximal domain of definition of its  diastasis $D^g_p$ at $p$, say $M_p$, is  dense on $M$. Simple examples of complete and   non-compact \K\ manifolds 
	where $M_p$ is not dense in $M$  can be found  in \cite[Theorem~1]{Suzuki82ContinuationDiastasis}.
	For instance, it is proven in  \cite{Suzuki82ContinuationDiastasis} that there exists a real analytic \K\ metric on $\C$ such that $M_0$ is a disk of finite radius. Notice also that one can exhibit examples of  real analytic \K\ metrics on $M$ where, for some $p\in M$, $D_p^g$ is negative at some point of its domain of definition. 
	The first example of such metrics was constructed by Calabi himself see   \cite[p.~23]{Calabi53Isometric}.
\end{ex}

\begin{remark}\rm\label{wbforever}
	We do not know any  example of  a   \K\ metric on a compact complex manifold which is not well-behaved at some point.  
	Moreover, in all known examples  with  $c_1(M)>0$ the maximal domain of definition of Calabi's diastasis function of a given \K\ metric at $p$
	is given by  the complement  of a complex analytic hypersurface $S_p$ of $M$ homologous to the dual of the \K\ form $\omega$ associated to $g$.
	Hence in these cases   condition $(A)$ in Definition~\ref{defWB} is always achieved by taking $H_p^g=S_p$.
	Notice that by Example~\ref{projindwb} this condition is satisfied  by all projectively induced  \K\ metrics regardless of the sign of $c_1(M)$.
	Calabi himself (in \cite[p.4]{Calabi53Isometric} 
	comments as: \lq\lq One could formulate several conjectures  over the behaviour of the diastasis in the large, which in the author's opinion would furnish ample material for further studies.''
\end{remark}

\begin{remark}\rm
	We also  lack examples  (compact or non-compact) where
	condition $(A)$ is valid  while   $(B)$ is violated.
\end{remark}

Further examples of well-behaved \K\ metrics 
satisfying conditions (i)-(iii) can be obtained by the following proposition on compactifications of $\C^n$ which will be a key ingredient in the proofs of Corollary~\ref{Cor1}-\ref{Cor3}.
By compactification of $\C^n$ we mean a compact complex manifold $M$ containing an analytic subvariety $H$ such that $X = M\setminus H$ is biholomorphic to $\C^n$.
We refer the reader to  \cite{PeternellSchneider91Compactification} for a survey on the topic.
\begin{prop}\label{proptoric}
	Let $M$ be a smooth projective compactification of $X$, such that $X$ is algebraically biholomorphic to $\C^n$ and let 
	$g$ be a \K\ real-analytic metric on $M$. Assume that there exists a point $p_*\in M$ where  the following conditions are satisfied:
	\begin{itemize}
		\item [$(A')$]
		Calabi's diastasis function $D^g_{p_*}$ around $p_*$ is globally defined and non-negative on $X$;
		\item [$(B')$]
		the Bochner coordinates  $(z)=\{z_1, \dots z_n\}$ around $p_*$ are globally defined on $X$.
	\end{itemize}
If  $\varphi: (S, h)\rightarrow (M, g)$ is a holomorphic isometric immersion from a compact \K\ manifold $(S, h)$ such that $\varphi (p)=p_*$, then $h$ is well-behaved and 
	satisfies (i)-(iii) of Theorem~\ref{mainteor}.
\end{prop}
\begin{proof}
	Condition $(A')$ and $(B')$ are clearly  stronger than $(A)$ and $(B)$ of Definition~\ref{defWB}. 
	Hence $g$ is well-behaved at $p_*$. 
	By the hereditary property of well-behaved \K\ metrics (cf. Example~\ref{hered}) $h$ turns out to be well-behaved at $p$. 
	Condition (i) is satisfied since $M$ is compact and condition (iii) is obviously hereditary.
	
	The least immediate part of the proof is to show the validity of (ii), namely the fact that the Bochner-Euclidean volume of $h$ at $p$ is infinite. 
	This follows by a modification of an argument given in the proof of \cite[Proposition 2.3]{ArezzoLoiZuddas16Toric}) valid for \K--Einstein submanifolds.
	Namely, let $w_1,\ldots,w_m$ be Bochner coordinates for $h$ defined in a neighbourhood $U$ of $p$.
	By Theorem~\ref{ThmCalabi}, we can assume that these are given as
		\begin{equation*}
		z_1\circ \varphi |_{U}=w_1,\dots , z_m\circ \varphi |_{U}=w_m.
		\end{equation*}
	Now consider the connected open subset $\hat{S}:=\varphi(S)\cap X$ of $\varphi(S)$. Denote by $\omega_g$ and $\omega_h$ the \K\ forms  associated to $g$ and $h$ respectively. 
	Now the two $m$-forms $\frac{\omega_g^m}{m!}$ and $\frac{i^m}{2^m}\det(g_{\alpha\bar\beta})\di z_1\wedge\di \bar z_1\wedge\cdots\wedge\di z_m\wedge\di \bar z_m$ defined on $X$ coincide on $\varphi(U)\subset\hat S$ as they both restrict to the volume form associated to $h$.
	Since they are real analytic, they must agree on the whole connected open set $\hat{S}$.
	
	Since $\frac{\omega_h^m}{m!}=\varphi^{*} \frac{\omega_g^m}{m!}$ is a volume form on $S\setminus H^h_p=\varphi^{-1}(\hat S)$, the form $\di z_1\wedge\di \bar z_1\wedge\cdots\wedge\di z_m\wedge\di \bar z_m$ cannot vanish on $\hat S$.
	We deduce that the restriction to $\hat{S}$ of the projection on the first $m$ coordinates $\pi:\C^n\lra\C^m$  is open. 
	Since it is also algebraic (because the biholomorphism between $X$ and $\C^n$ is algebraic), its image $\pi(\hat{S})$ contains a Zariski open subset of $\C^m$, see \cite[Theorem~13.2]{Borel91LinearGroups}.
	We conclude that the Bochner--Euclidean volume is infinite because
	$$  \vol_{eucl} (S\setminus H^h_p, h)\geq\int_{\pi(\hat{S})}\dfrac{i^m}{2^m}\di z_1\wedge\di \bar z_1\wedge\cdots\wedge\di z_m\wedge\di \bar z_m=\infty.$$
\end{proof}

\begin{remark}\label{RmkEsempioCinCPn}\rm
	The compactness assumption of $S$ in Proposition~\ref{proptoric} is necessary since (i) of Theorem~\ref{mainteor}, unlike (ii) and (iii), is not hereditary for complete \K\ manifold
	admitting a holomorphic isometry into
	finite dimensional complex projective space. 
	Take, for example, a complex $2$-dimensional abelian variety $T\subset \C P^n$ and choose a copy of $	\C$ dense in $T$. Restrict the Fubini-Study metric of  $\C P^n$
	to $T$ and to  $\C$ denoting the later by  $g$. Thus $(\C, g)$ is a non-compact  and complete
	\K\  manifold of infinite volume,  which admits a holomorphic isometric immersion into $(\C P^n, g_{FS})$ (see \cite[Example 6.4, p. 18]{LoiPlaciniZedda23SasakiHomogeneous} for details).
\end{remark}

\begin{remark}\rm
	It is worth pointing out that property
	(ii) is not  hereditary when the ambient space is not a compactification of $\C^n$.
	Take, for example two complex manifolds $M_1$ and $M_2$ equipped with two \K\ metrics $g_1$ and $g_2$ which are well-behaved at $p_1$ and $p_2$ respectively. Choose $(M_1,g_1)$ and $(M_2,g_2)$ satisfying $\vol_{eucl}(M_1\setminus H_{p_1}^{g_1})<\infty$ and $\vol_{eucl}(M_2\setminus H_{p_2}^{g_2})=\infty$.
	Then the product $M_1\times M_2$ equipped with the \K\ metric $g_1\oplus g_2$  satisfies (iii) while  $g_1$ (which is induced by the natural inclusion $M_1\rightarrow M_1\times M_2$)  does not.
\end{remark}

The prototype example of compactifications of $\C^n$ satisfying 
$(A')$ and $(B')$ is  of course 
the finite dimensional complex projective space $\C P^n$.
Indeed, we can write $\C P^n=\C^n\sqcup H^{g_{FS}}_{p_0}$ (cf. Example~\ref{projspace}). 
Hence, combining Example~\ref{projindwb} with  Proposition~\ref{proptoric} proves Corollary~\ref{Cor1}.

Notice that by Proposition~\ref{proptoric}  all   the \K\ metrics $g$ in Example~\ref{genflag} and Example~\ref{toric} respectively  and all their \K\ submanifolds are (well-behaved) at a suitable $p$ and satisfy (i)-(iii)
of Theorem~\ref{mainteor}. Namely, we have proven Corollary~\ref{Cor2} and Corollary~\ref{Cor3}.

\section{Proof of Theorem~\ref{mainteor} and a toy example}\label{SecProofs}
\begin{proof}[Proof of Theorem~\ref{mainteor}]
Take Bochner coordinates $\{z_1,\dots ,z_n\}$ for the metric $g$ on a contractible  neighbourhood $U$ of $p$.
Then   
$$\omega =\frac{i}{2}\sum g_{\alpha\bar{\beta}}\di z_{\alpha}\wedge \di\bar{z}_{\beta}=\frac{i}{2}\sum\frac{\partial ^2 D^g_p}
{\partial z_{\alpha} \partial \bar z_{\beta}} \di z_{\alpha}\wedge \di\bar{z}_{\beta}$$ on $U$
and 
$\rho _{\omega}=-i\partial \bar{\partial}\log \det g_{{\alpha}\bar{\beta}}$ is the local expression of its \emph{Ricci form}.
From the $\partial\bar\partial$-lemma and the equation $\rho_\omega=\lambda\Omega$ it follows that the volume form of $(M, g)$ reads on $U$ as:
\begin{equation}\label{voleuclbo}
\frac{\omega^n}{n!}=
\det \left(\frac{\partial ^2 D^g_p}
{\partial z_{\alpha} \partial \bar z_{\beta}}\right)
d\mu(z)=
e^{-\frac{\lambda}{2}D^G_p+F+\bar F}
d\mu(z)\, \, ,
\end{equation}
where  $d\mu(z):=\frac{i^n}{2^n}dz_1\wedge d\bar z_1\wedge\dots\wedge 
dz_n\wedge d\bar z_n$, 
 $F$ is a holomorphic
function on $U$
and $D_p^g$ 
(resp. $D^G_p$)
is the diastasis function with respect to $p$
for the \K\ metric $g$ (resp. $G$).
 
We claim that
$F+\bar F=0$.
In order to prove our  claim notice 
that, by  the very definition
of Bochner coordinates, one can easily
check that 
$\log\det (\frac{\partial ^2 D^g_p}
{\partial z_{\alpha} \partial \bar z_{\beta}})$
is of diastasis type.
Moreover,  formula  \eqref{voleucl} 
yields
$$F + \bar F=
\frac{\lambda}{2} D^G_p+
\log\det (\frac{\partial ^2 D^g_p}
{\partial z_{\alpha} \partial \bar z_{\beta}})$$
where the  right hand side is 
of diastasis type because $D^G_p$ is of diastasis type
by Remark~\ref{remdefdiast}. This forces
$F + \bar F$ to vanish,
proving our claim.

By assumption (B) of Definition~\ref{defWB}  there exist  $f_1, \dots ,f_n$  holomorphic functions on $M\setminus H_p^g$ extending 
the Bochner coordinates $\{z_1, \dots ,z_n\}$ for the metric $g$.
On the other hand, since assumption $(A)$ of Definition~\ref{defWB}    is satisfied both for $g$ and $G$,
the real analytic  $n$-forms 
$\frac{\omega^n}{n!}$
and
$e^{-\frac{\lambda}{2}D^G_{p}}
d\mu (f)$ 
are globally defined on the connected   open set 
$$X:=(M\setminus H_p^g)\cap (M\setminus H_p^G)=M\setminus (H_p^g\cup H_p^G).$$
Notice that $X$ is connected by assumptions (iii) and (iv).
Since, by 
formula \eqref{voleucl} with $F+\bar F=0$, these $n$-forms agree on $U$,
they
must agree
on $X$, i.e.
\begin{equation}\label{voleucl3}
\frac{\omega^n}{n!}=
e^{-\frac{\lambda}{2}D^G_p}
d\mu (f).
\end{equation}

Now assume by contradiction that $\lambda\leq 0$ .
Combining
\eqref{voleucl3} and the fact that $D^G_{p}$ is non-negative on $X$
(by assumption $(A)$ for the metric $G$) 
yields
$$\vol (M, g):=\int_M\frac{\omega^n}{n!}=\int_X\frac{\omega^n}{n!}=\int_Xe^{-\frac{\lambda}{2}D^G_p}
d\mu (f)\geq\int_X d\mu (f).$$
Moreover, since $g$ has infinite Bochner-Euclidean volume (assumption (ii)), we get

$$\vol (M, g)\geq\int_X d\mu (f)=\int_{M\setminus H_p^g} d\mu (f)=\vol_{eucl} (M\setminus H_p^g, g)=\infty,$$
which is the desired contradiction,
being $\vol (M, g)$ finite,  by assumption (i).
\end{proof}

\begin{remark}\rm\label{rmkafterth}
Notice that in the proof of Theorem~\ref{mainteor} we do not use the non-negativity of $D_p^g$  of condition 
$(A)$ for the  metric $g$,
but only the fact that $D^g_p$ is globally defined on $M\setminus H_p^g$. 
Moreover, in the proof  we are  only requiring $D^G_p$ to satisfy condition $(A)$ and not necessarily $(B)$ of  Definition~\ref{defWB}.
Therefore we get:
\begin{cor}\label{CorConditionsProjInduced}
	Let $g$ be a projectively induced metric on a compact \K\ manifold $M$ with $c_1(M)<0$. If $-\rho_\omega$ is associated to a \K\ metric $G$, then either $G$ does not satisfy $(A)$ of Definition 
	~\ref{defWB}  or 
	$\codim_{\R}(H_p^G)\leq 1$.
\end{cor}
\end{remark}
As a first example we specialize to \K-Ricci iterations on simply connected compact homogeneous \K\ manifolds, that is, generalized flag manifolds.
This choice is based on the fact that homogeneous \K\ manifolds are a precious source of examples of \K\ manifolds with positive first Chern class in the compact realm.
In particular, to our best knowledge, the only examples of projectively induced \K\ forms $\omega$, $\Omega$ satisfying equation~\eqref{neweq} and the hypotheses of Theorem~\ref{mainteor} are homogeneous \K\ forms.
In fact, we can say something more for generalized flag manifolds. 

\begin{prop}\label{PropHomogeneous}
	Let $M=G/K$ be a generalized flag manifold endowed with a \K\ metric $g$ and associated \K\ form $\omega$.
	Then $\omega$ is homogeneous if and only if its Ricci form $\rho_\omega$ is a \K-Einstein form.
\end{prop}
\begin{proof}
	Let $\omega$ be homogeneous, i.e., for any $h\in G$ we have $h^*\omega=\omega$.
	Then we have $h^*\rho_\omega=\rho_{h^*\omega}=\rho_\omega$ so that $\Omega:=\rho_\omega$ is a homogeneous form.
	Moreover, $\Omega$ is positive at a point because it represents a positive class.
	This, together with homogeneity, shows that $\Omega$ is a \K\ form.
	Thus  $\Omega$ \K-Einstein because it is has constant scalar curvature and satisfies $[\rho_\Omega]=[\Omega]$, see for instance \cite[Section~2]{Kobayashi67JDGcscK}.
	
	Vice versa, suppose that $\omega$ is any \K\ form on $M$ such that $\rho_\omega$ is homogeneous.
	Then we have $\rho_{h^*\omega}=h^*\rho_\omega=\rho_\omega$.
	Since $G$ is connected, $h^*=\Id:H^*(M)\lra H^*(M)$ because $h$ is homotopic to the identity in $G$ as a map $M\lra M$.
	This in particular implies $[h^*\omega]=h^*[\omega]=[\omega]$ and we can conclude that $\omega=h^*\omega$ by Calabi's conjecture.
\end{proof}

	Observe that an integral homogeneous \K\ metric on a generalized flag manifold is always projectively induced, see for instance \cite{LoiZuddasGromovWidthHKM,Takeuchi78Homogeneous}.
	It seems natural to ask whether, relaxing the homogeneity condition on the \K\ metric, the statement of Proposition~\ref{PropHomogeneous} still holds true. Namely, we ask the following question:
	\vskip 0.2 cm
	\noindent \textbf{Question:} \emph{Let $g$ and $G$ be projectively induced \K\ metrics on a generalized flag manifold. Is $\rho_\omega=\lambda\Omega$ equivalent to $G$ being \K--Einstein?}
		\vskip 0.2 cm
	As a special case, Proposition~\ref{PropProjIndCP1} below  gives an answer in the affirmative for  (radial) metrics on $\CP^1$. In order to prove this result, we carry out in full detail the K\"ahler-Ricci iterations on projectively induced metrics on $\CP^1$.
Namely, suppose that a \K\ metric $g$ is induced by a holomorphic embedding $\varphi:\CP^1\lra\CP^n$ in a higher dimensional complex projective space, that is, $g=\varphi^*g_{FS}$ where $g_{FS}$ is the Fubini-Study metric of $\CP^n$.

For homogeneous coordinates $[z_0:z_1]$ consider the affine chart $\{z_0\neq0\}$ with coordinate $z=\frac{z_1}{z_0}$. 
Recall that the embedding $\varphi$ is determined by sections of a bundle $\OO(k)$ which can be identified with homogeneous polynomials of degree $k$ in $z_0$ and $z_1$. 
Without loss of generality we can assume that $\varphi$ is full, i.e., the embedding is determined by a basis of $H^0(\CP^1,\OO(n))$.
Moreover, let us assume that the induced metric $\varphi^*(g_{FS})$ is radial in  $z$.
This is equivalent to the fact that the coordinates of $\varphi$ are monomials so that in homogeneous coordinates (possibly after reordering) $\varphi$ is of the form
$$\varphi([z_0:z_1])=[\alpha_0z_0^n:\alpha_1z_1z_0^{n-1}:\cdots:\alpha_nz_1^n].$$
Notice that we can impose $\alpha_0=\alpha_n=1$ after precomposing $\varphi$ with $f\in PGL(2,\C)$ given in coordinates by $f([z_0:z_1])=\left[\frac{1}{\sqrt[n]{\alpha_0}}z_0:\frac{1}{\sqrt[n]{\alpha_n}}z_1\right]$.
Then the \K\ form $\omega=\varphi^*(\omega_{FS})$ associated to the projectively induced metric $g$ is given in the chosen chart by
\begin{equation}\label{EqMetricProjInduced}
	\omega=\dfrac{i}{2}\de\deb\log(1+a_1\vert z\vert^2+\cdots+a_j\vert z\vert^{2j}+\cdots+\vert z\vert^{2n})
\end{equation}
where $a_j=\vert\alpha_j\vert^2$.
Denote the argument of the logarithm in \eqref{EqMetricProjInduced} by $P\left(x\right)$ where $x=\vert z\vert^2$.
Then the form can be rewritten as 
\begin{equation}\label{EqMetricProjInduced2}
	\omega=\dfrac{i}{2}\dfrac{P'P+xP''P-x(P')^2}{P^2}\di z\wedge\di \bar z.
\end{equation}
Writing $D:=P'P+xP''P-x(P')^2$ for conciseness, the Ricci form $\rho_\omega $ is given by
\begin{equation}\label{EqRicci1}
	\rho_\omega =-i\de\deb\log\dfrac{D}{2P^2}=i\de\deb\log\dfrac{P^2}{D}.
\end{equation}
First we want to prove the following result.
\begin{prop}\label{PropProjIndCP1}
	Consider a projectively induced metric $g$ on $\CP^1$ with associated \K\ form $\omega$ as in \eqref{EqMetricProjInduced}. Then its Ricci form $\rho_\omega  $ is projectively induced if and only if $a_j={{n}\choose {j}}$, i.e.,  $g=ng_{FS}.$ 
\end{prop}
In order to prove Proposition~\ref{PropProjIndCP1} we need the following auxiliary result.
\begin{lem}\label{LemProjIndCP1}
	Let $\omega=\frac{i}{2}\de\deb\log\dfrac{f(x)}{h(x)}$ be a \K\ form on $\CP^1$ where $z$ is the affine coordinate on the chart $U_0=\{z_0\neq0\}$, $x=\vert z\vert^2$ and $f,h\in\R[x]$ are polynomials.
	If $g$ is projectively induced, then $f/h=Q\in\R[x]$.
\end{lem}
\begin{proof}
	Suppose there exists an embedding $\varphi:\CP^1\lra\CP^n$ such that $\omega=\varphi^*\omega_{FS}$.
	Since $\omega$ is a radial form on $U_0$, its potential must be radial on $U_0$. This means that we can write $$\omega=\frac{i}{2}\de\deb\log Q(x)= \frac{i}{2}\de\deb\log\dfrac{f(x)}{h(x)}$$
	where $\varphi(z)=[\varphi_0(z):\cdots:\varphi_n(z)]$ and $Q=\sum\vert\varphi_j(z)\vert^2$.
	Thus, $\log(f/h)$ differs from $\log Q$ by the real part of a holomorphic function $F:U_0\lra\C$, that is, $Q=f/h+e^{F+\bar F}$.
	Now it is clear from the series expansion of $Q$ that $F+\bar F$ must be constant and this proves the claim.
\end{proof}

\begin{proof}[Proof of Proposition~\ref{PropProjIndCP1}]
	Clearly, if $P(x)=\left(1+x\right)^n=\left(1+\vert z\vert^2\right)^n$ then $\rho_\omega =\rho(n\omega_{FS})=4\omega_{FS}$ which is again projectively induced.
	
	Vice versa, suppose $\rho_\omega = i\de\deb\log\dfrac{P^2}{D}$ is a projectively induced \K\ form.
	Then by Lemma~\ref{LemProjIndCP1} we can assume that 
	\begin{equation}\label{EqFraction}
		\dfrac{P^2}{D}=\lambda Q
	\end{equation}
	for a monic polynomial $Q$ and $\lambda\in\R$. 
	Notice that $\frac{i}{2\pi}\rho_\omega $ represents $c_1(\CP^1)$. Therefore $\deg Q=2$ for topological reasons.
	
	We want to show that $P=(1+x)^n$ by contradiction. 
	If $P\neq(1+x)^n$, then $P$ must have more than one root. 
	Consequently there exists a root of $P$, say $\alpha$, such that $Q\neq(x-\alpha)^2$.
	To conclude the proof we show that this implies $P=(x-\alpha)^n$ which is the desired contradiction.

	We do so by proving that if $(x-\alpha)^{k}$ divides $P$ for $k<n$, then also $(x-\alpha)^{k+1}$ divides $P$.
	First, it is evident from \eqref{EqFraction} that if $(x-\alpha)^{k}$ divides $P$, then $(x-\alpha)^{2k-1}$ divides $D$ because $\alpha$ is at most a simple root of $Q$.
	Now it is easy to see that the $(2k-2)$-th derivative $D^{(2k-2)}$ of $D$ can be written as
	\begin{equation}\label{EqDerivativesPoly}
		D^{(2k-2)}=-x\left(P^{(k)}\right)^2+R
	\end{equation}
	where $P^{(k)}$ is the  $k$-th derivative of $P$ and $R$ is a sum of polynomials each of which is divided by a derivative of $P$ of order $< k$. 
	Therefore we get
	\begin{equation}\label{EqDerivativesPolyEvaluated}
		0=D^{(2k-2)}(\alpha)=-\alpha\left(P^{(k)}(\alpha)\right)^2+R(\alpha)=-\alpha\left(P^{(k)}(\alpha)\right)^2
	\end{equation}
	Since $P=(1+\cdots)$ we have $\alpha\neq0$ and $P^{k}(\alpha)=0$, that is, $(x-\alpha)^{k+1}$ divides $P$.
	This concludes the proof.
\end{proof}

We want to specialize now to the case of embeddings $\varphi:\CP^1\lra\CP^2$ in order to give some explicit examples. 
\begin{ex}\label{ExPIntoWB}\rm
	By the discussion above, a radial embedding $\varphi:\CP^1\lra\CP^2$ can be written in homogeneous coordinates in the form
	$$\varphi([z_0:z_1])=[z_0^2:\alpha z_1z_0:z_1^2].$$
	The \K\ form $\omega$ associated to the projectively induced metric $g=\varphi^*(g_{FS})$ metric on $\CP^1$ in an affine chart is given as
	\begin{equation}\label{EqMetricProjInducedCP2}
		\omega=\dfrac{i}{2}\de\deb\log(1+ax+x^2)=\dfrac{i}{2}\dfrac{(a+4x+ax^2)}{(1+ax+x^2)^2}\di z\wedge\di\bar{z}
	\end{equation}
	where $a=\vert\alpha\vert^2$ and $x=\vert z\vert^2$. 
	Therefore its Ricci form is
	\begin{equation}\label{EqMetricWB1Explicit}
		\omega_1=\rho_\omega=-i\de\deb\log\dfrac{(a+4x+ax^2)}{2(1+ax+x^2)^2}=i\dfrac{2(a+4x+ax^2)^3-4a(1+ax+x^2)^3}{(a+4x+ax^2)^2(1+ax+x^2)^2}\di z\wedge\di\bar{z}.
	\end{equation}
	By Proposition~\ref{PropProjIndCP1} the metric $g_1$ associated to the Ricci form $\omega_1:=\rho_\omega$ is projectively induced if and only if $a=2$. 
	Moreover, one can easily check that the coefficients of the polynomial $2(a+4x+ax^2)^3-4a(1+ax+x^2)^3$ are all positive for $\sqrt{2}\leq a\leq2$ so that for those values $g_1$ is indeed a \K\ metric.
	In addition, the expression \eqref{EqMetricWB1Explicit} of the \K\ form $\omega_1$ is written in Bochner coordinates and can be used to compute the Euclidean volume $\vol_{eucl}(\CP^1\setminus [0:1], g_1)
	$:
	
	$$\vol_{eucl}(\CP^1\setminus [0:1], g_1)=2\pi\int_{0}^\infty\dfrac{2(a+4x+ax^2)^3-4a(1+ax+x^2)^3}{(a+4x+ax^2)^2(1+ax+x^2)^2}\di x$$
	where we made use of the equality $\di x=\di r^2=2r\di r$ for $z=re^{-i\theta}$. Looking at the degrees of the numerator and denominator of the integrand it is immediate that $\vol_{eucl}(\CP^1\setminus [0:1], g_1)=\infty$.
	We conclude that $g$ is a projectively induced metric, hence well-behaved at all points and with infinite euclidean volume, whose Ricci form is associated to a well-behaved, not projectively induced \K\ metric $g_1$ with infinite Euclidean volume for $\sqrt{2}\leq a<2$.
	
	Based on this example we want to show that the second Ricci iteration of $g$ is again a \K\ form $\omega_2$ associated to a well-behaved \K\ metric $g_2$ with infinite Euclidean volume when $a$ belongs to an suitable interval.
	To show this we iterate the computations carried above. 
	It is evident that the complexity of such calculations increases with the number of iterations. 
	For this reason we resorted to the use of the software Mathematica \cite{Mathematica} to compute the Ricci form $\omega_2:=\rho_{\omega_1}=\rho^2_\omega$ as well as the values of $a$ for which $\omega_2$ is associated to a \K\ metric.
	The calculations yield
	\begin{align}\label{EqMetricWB2Explicit}
		\nonumber	\omega_2=\rho^2_\omega=-i\de\deb\log\dfrac{2(a+4x+ax^2)^3-4a(1+ax+x^2)^3}{(a+4x+ax^2)^2(1+ax+x^2)^2}=i\dfrac{6N}{D}\di z\wedge\di\bar{z}
	\end{align}
	where $N$ and $D$ are polynomials in $x$ with $D>0$. 
	One can check numerically that the coefficients of $N$ are all positive when $a$ lies in an interval $I\subset(\sqrt{2},2)$.
	Therefore, for $a\in I$, both $g_1$ and $g_2$ are positive definite.
	Moreover, from $\deg N=18$ and $\deg D=20$ follows that the Euclidean volume of $g_2$ is
	$$\vol_{eucl}(\CP^1\setminus [0:1], g_2)=2\pi\int_{0}^\infty\dfrac{6N}{D}\di x=\infty.$$
	Therefore, for $a\in I$, $g_1$ is an example of well-behaved, not projectively induced \K\ metric with infinite Euclidean volume whose Ricci form is associated to a  well-behaved, not projectively induced \K\ metric with infinite Euclidean volume. 
	This provides a pair of well-behaved (not projectively induced) metrics on $\CP^1$ satisfying the hypotheses of Theorem~\ref{mainteor}. 
\end{ex}

\section{Proofs of Theorem~\ref{mainteor2} and Theorem~\ref{mainteor3}}\label{SecFinal}

In order to prove Theorem~\ref{mainteor2} we briefly recall some basic facts on Umehara algebra and its field of fractions.
Let $M$ be a complex manifold. Fix a point  $p\in M$ and let 
$\OO_p$ be the algebra of germs of holomorphic functions around $p$.
Denote  by $\R_p$ the germs of real numbers. 
The Umehara algebra (see  \cite{LoiMossa21KRSIntoComplexSpaceForms,Umehara88Diastases} and references therein) is defined  to be the $\R$-algebra $\Lambda_p$ generated
by elements of the form  $h\bar k+\bar h k$, for $h, k\in \OO_p$. 
Let
$$\hat\OO_p=\left\{\alpha =(\alpha_1,   \dots ,\alpha_m) \ | \   \alpha_j\in\OO_p,\ \alpha_j(p)=0 , \forall j=1, \dots , m, m\geq 1\right\}.$$
For  $\alpha =(\alpha_1,  \dots ,\alpha_m)\in \hat\OO_p$
and $\ell\in\N$ such that $\ell\leq m$
we set
$$\langle \alpha, \alpha\rangle_{\ell}(z)=\sum_{j=1}^\ell|\alpha_j(z)|^2- \sum_{k=\ell+1}^{m}|\alpha_k(z)|^2.$$
Since $h\bar k+\bar h k=|h+k|^2-|h|^2-|k|^2$, it is not hard to see 
that each $f \in \Lambda_p$  can be written  as
$$f=h+\bar h+\langle \alpha, \alpha\rangle_{\ell}$$
for some $h\in \OO_p$, $\alpha =(\alpha_1,  \dots ,\alpha_m)\in \hat\OO_p$, 
$\ell\leq m$
and such that  $\alpha_1, \dots , \alpha_m$  are linearly independent over $\C$.

Consider  the $\R$-algebra $\tilde\Lambda_p\subset \Lambda_p$ given by
\begin{equation}\label{Lambda0}
\tilde\Lambda_p=\left\{a+ \langle \alpha, \alpha\rangle_{\ell} \mid a\in \R_p,\ \alpha\in \hat\OO_p,\  \ell\leq m\right\}.
\end{equation}
Notice that the germ of the real part of a non-constant  holomorphic function 
$h\in\OO_p$ belongs to 
$\Lambda_p$ but not to 
$\tilde\Lambda_p$. 

The key elements in the proof of Theorem~\ref{mainteor2} and Theorem~\ref{mainteor3} are  the following two  lemmata.
 The first lemma descents from  \cite[Theorem~2.1]{LoiMossa22HBD} and we specialize it here to our needs.

\begin{lem}\label{LemLogConst}
Let $\tilde K_p$ be the field of fractions of  $\tilde\Lambda_p$.
Let $\xi\in \tilde K_p$, then
$e^\xi\not\in \tilde K_p^\alpha\tilde K_p\setminus \R_p ,$
$\forall \alpha\in\R$.
\end{lem}

\begin{lem}\label{umeharaimp}
Let $M$ be  an $n$-dimensional complex manifold and  $p\in M$ 
and let $\frac{f}{h}\in K_p$, where $K_p$ denotes the field of fractions of the Umehara algebra 
$\Lambda_p$. Then,  for any system of complex coordinates $\{z_1, \dots ,z_n\}$ around $p$ one has:
$$f^{n+1}h^{n+1}\det\left[\frac{\partial^2\log\frac{f}{h}}{\partial z_{\alpha}\partial\bar z_{\beta}}\right]\in \Lambda_p, \ \mbox{with } \alpha, \beta =1, \dots , n.$$
\end{lem}
\begin{proof}
Set  $f_\alpha=\partial f / \partial z_\alpha$, $h_\alpha=\partial h / \partial z_\alpha$, $f_{\bar{\beta}}=\partial f / \partial \bar{z}_\beta$, $h_{\bar{\beta}}=\partial h / \partial \bar{z}_\beta$ and 
$f_{\alpha \bar{\beta}}=\partial^2 f / \partial z_\alpha \partial \bar{z}_\beta$, $h_{\alpha \bar{\beta}}=\partial^2 h / \partial z_\alpha \partial \bar{z}_\beta$, $\alpha, \beta =1, \cdots, n$. Then
\begin{equation*}
\frac{\partial^2 \log(f/h)}{\partial z_\alpha \partial \bar{z}_\beta}
=  f_{\alpha \bar{\beta}} / f
- h_{\alpha \bar{\beta}} / h
-f_\alpha f_{\bar{\beta}} / f^2
+h_\alpha f_{\bar{\beta}} / h^2,
\end{equation*}
so that 
\begin{equation}\label{EqMatrixEntry}
fh	\frac{\partial^2 \log(f/h)}{\partial z_\alpha \partial \bar{z}_\beta}
	= \left(h\left(f_{\alpha \bar{\beta}} -f_\alpha f_{\bar{\beta}} / f\right)  -f\left(h_{\alpha \bar{\beta}}- h_\alpha h_{\bar{\beta}} / h\right)\right) .
\end{equation}
Now some linear algebra shows that  $f^{n+1}h^{n+1}\det\left[\frac{\partial^2\log\frac{f}{h}}{\partial z_{\alpha}\partial\bar z_{\beta}}\right]$ is the determinant of a matrix whose entries are generated by holomorphic and anti-holomorphic functions.
Namely, we have
\begin{equation*}
f^{n+1}h^{n+1} \operatorname{det}\frac{\partial^2 \log(f/h)}{\partial z_\alpha \partial \bar{z}_\beta}
=fh \operatorname{det}\left(\begin{array}{ccccc}
	fh\frac{\partial^2 \log(f/h)}{\partial z_1 \partial \bar{z}_1}& \cdots & fh\frac{\partial^2 \log(f/h)}{\partial z_1 \partial \bar{z}_n}& 0 & 0\\
	\vdots & \ddots & \vdots& \vdots & \vdots \\
	fh\frac{\partial^2 \log(f/h)}{\partial z_n \partial \bar{z}_1}& \cdots & fh\frac{\partial^2 \log(f/h)}{\partial z_n \partial \bar{z}_n}& 0 & 0\\
	hf_{\bar{\mathrm{1}}} / f & \cdots & h f_{\bar{n}} / f & 1& 0\\
	fh_{\bar{\mathrm{1}}} / h & \cdots & f h_{\bar{n}} / h & 0& 1
\end{array}\right) .
\end{equation*}
Denote by $R_j$ the $j$-th row of the matrix above. Replacing $R_j$ by $R_j+f_jR_{n+1}-h_jR_{n+2}$ for all $j=1,\ldots,n$ and making use of \eqref{EqMatrixEntry} we get
\begin{equation*}
	f^{n+1}h^{n+1} \operatorname{det}\frac{\partial^2 \log(f/h)}{\partial z_\alpha \partial \bar{z}_\beta}
=fh \operatorname{det}\left(\begin{array}{ccccc}
	f_{1 \overline{1}}h-fh_{1 \overline{1}} & \cdots & f_{1 \bar{n}}h-fh_{1 \overline{n}} & f_1& -h_1 \\
	\vdots & \ddots & \vdots & \vdots & \vdots \\
	f_{n \overline{1}}h-fh_{n \overline{1}} & \cdots & f_{n \bar{n}}h-fh_{n \overline{n}} & f_n & -h_n \\
	hf_{\overline{1}} / f & \cdots &h f_{\bar{n}} / f & 1 & 0 \\
	fh_{\overline{1}} / h & \cdots & fh_{\bar{n}} / h & 0 & 1
\end{array}\right).
\end{equation*}
Finally, multiplying the last two rows by $f$ and $h$ respectively we get
\begin{equation*}
	f^{n+1}h^{n+1} \operatorname{det}\frac{\partial^2 \log(f/h)}{\partial z_\alpha \partial \bar{z}_\beta}
	= \operatorname{det}\left(\begin{array}{ccccc}
		f_{1 \overline{1}}h-fh_{1 \overline{1}} & \cdots & f_{1 \bar{n}}h-fh_{1 \overline{n}} & f_1& -h_1 \\
		\vdots & \ddots & \vdots & \vdots & \vdots \\
		f_{n \overline{1}}h-fh_{n \overline{1}} & \cdots & f_{n \bar{n}}h-fh_{n \overline{n}} & f_n & -h_n \\
		hf_{\overline{1}} & \cdots &h f_{\bar{n}}& f & 0 \\
		fh_{\overline{1}} & \cdots & fh_{\bar{n}}  & 0 & h
	\end{array}\right).
\end{equation*}
Hence $f^{n+1}h^{n+1} \operatorname{det}\frac{\partial^2 \log(f/h)}{\partial z^\alpha \partial \bar{z}^\beta}$ is finitely generated by holomorphic or anti-holomorphic functions. In addition, it is real-valued, because the matrix $\frac{\partial^2 \log(f/h)}{\partial z^\alpha \partial \bar{z}^\beta}$ is Hermitian and this proves the lemma.
\end{proof}
\begin{remark}\rm
Lemma~\ref{umeharaimp} is an extension  of  \cite[Lemma 2.2]{Umehara87Einstein}, which is valid for $\Lambda_p$, to its field of fractions $K_p$.
\end{remark}

The following result, needed both in the proofs of Theorem~\ref{mainteor2}
and Theorem~\ref{mainteor3}, follows  by induction on $k\geq 1$,  by the definition of the Ricci form associated to a real analytic \K\ metric and by Lemma~\ref{umeharaimp}.

\begin{cor}\label{corrhok}
Let $g$ be a real analytic \K\ metric on a complex manifold $M$.
Choose $U$ around $p\in M$ such that 
Calabi's diastasis function $D_p^g:U\rightarrow \R$ is defined 
and such that 
$\rho^k_\omega=\frac{i}{2}\partial\bar\partial \Psi_k$, where 
$\Psi_k$ is of diastasis type (cf. Remark~\ref{locexp(1,1)}).
If $D_p^g$ belongs either to
$\tilde K_p$ or  $\log\tilde K_p$ then $\Psi_k\in\log\tilde K_p$.
\end{cor}

\begin{proof}[Proof of Theorem~\ref{mainteor2}]
Since $g$ is induced by a flat metric there exists 
an open neighbourhood $U$ of a point $p\in M$ and a holomorphic isometry $\varphi:U\rightarrow \C^N$ such that $\varphi^*g_0=g_{\vert_U}$,
where $g_0$ is the flat metric on $\C^N$.
Let $$\varphi_{|U}:U\rightarrow \C^N, q\mapsto (\varphi_1(q), \dots , \varphi_N(q))$$
where $\varphi_j\in \CO_p$ and $\varphi_j(p)=0$, $j=1, \dots ,N$,
be the local expression of $\varphi$.

Notice that, by the hereditary property of the diastasis function  we  have 
\begin{equation}\label{Dpflat}
D^g_p =\sum_{i=1}^N \left|\varphi_i\right|^2\in\tilde\Lambda_p.
\end{equation}
By shrinking $U$ if necessary,  Corollary~\ref{corrhok} yields
\begin{equation}\label{ricciformkind}
\rho^k_\omega=\frac{i}{2}\partial\bar\partial\Psi_k
\end{equation} 
where $\omega$ is the \K\ form associated to $g$
and $\Psi_k\in\log\tilde K_p$.
Now equation $\rho^k_\omega=\lambda\Omega$, together 
with the fact that  $G$ is induced by the flat metric
gives
$$\Psi_k= \lambda D_p^G=\lambda\sum_{i=1}^M \left|\phi_i\right|^2\in\tilde\Lambda_p,$$
 where $\phi_k\in \CO_p$ and $\phi_k(p)=0$, 
 such that $\phi_k$ is  a  non-constant function for all $k=1, \dots ,M$.
 Thus Lemma~\ref{LemLogConst} forces  $\lambda =0$.
 Then the Ricci tensor of $g$ vanishes.
 Since the immersion of $M$ in its flat ambient space lifts to a globally defined holomorphic isometry $\varphi :M\rightarrow \C^N$, the thesis follows.
 \end{proof}
 
 \begin{remark}\rm
 In contrast with the proof of Theorem~\ref{mainteor3}, the previous proof only uses a weak version of Lemma~\ref{LemLogConst}, that is, the fact that 
$e^f\not\in\tilde K_p\setminus \R_p $, 
$\forall f \in\tilde\Lambda_p$. 
This has been proven in \cite[(i) of Theorem~2.1)]{ChengDiscalaYuan17UmeharaAlgebra}.
\end{remark}

\begin{proof}[Proof of Theorem~\ref{mainteor3}]
Fix a point $p\in M$ and consider an open coordinate neighbourhood centred at $p$..
Since the solitonic vector field $X$ for $G$ can be assumed to be the real part of a holomorphic vector field, we can write locally
$$
X=\sum_{j=1}^{n}\left(f_{j} \frac{\partial}{\partial z_j}+\bar{f}_{j}\frac{\partial}{\partial\bar z_j}\right)
$$
for some holomorphic functions $f_{j}, j=1, \ldots, n$, on $U$.
Thus, by the definition of Lie derivative, after a 
straightforward computation, we can write on $U$
\begin{equation}\label{LXomega}
L_X\Omega=\frac{i}{2}\partial\bar\partial\xi
\end{equation}
where 
\begin{equation}\label{fXlocal}
\xi=\sum_{j=1}^nf_j\frac{\partial D^G_p}{\partial z_j}+\bar f_j\frac{\partial D^G_p}{\partial \bar z_j}.
\end{equation}
The equation
$\rho^k_\omega=\lambda\Omega$
together with the KRS  equation for $\Omega$, i.e. $\rho_\Omega =\mu\Omega+L_X\Omega$,  gives
\begin{equation}\label{solitonk}
\rho^{k+1}_{\omega}=\mu\rho^k_\omega +L_X\Omega,
\end{equation}
By Remark~\ref{locexp(1,1)} we can restrict $U$ so that $\rho^k_\omega=\frac{i}{2}\partial\bar\partial\Psi_k$ and 
$\rho^{k+1}_\omega=\frac{i}{2}\partial\bar\partial\Psi_{k+1}$ with $\Psi_k$ and $\Psi_{k+1}$
of diastasis type.
Thus \eqref{LXomega}, \eqref{fXlocal}  and \eqref{solitonk} together with $\lambda D_p^G=\Psi_k$ yield
\begin{equation}\label{Eqfinale}
\Psi_{k+1}=\lambda\mu\Psi_k+\dfrac{\xi}{\lambda}.
\end{equation}
Notice now that the assumption that $g$ is induced by a complex space form implies that $D_p^g\in \tilde K_p$ (in the flat case)
or $D_p^g\in \log\tilde K_p$ (for the hyperbolic or projective space). Thus, by Corollary~\ref{corrhok} $\Psi_k$ and $\Psi_{k+1}$
belong to $\log\tilde K_p$ and, consequently,  $\dfrac{\xi}{\lambda}=\sum_{j=1}^nf_j\frac{\partial \Psi_k}{\partial z_j}+\bar f_j\frac{\partial\Psi_k}{\partial \bar z_j}\in\tilde K_p$.
It follows by \eqref{Eqfinale} that
$e^{\xi}\in \tilde K_p^{-\lambda\mu}\tilde K_p$. 
Hence $\xi$ is a constant by Lemma~\ref{LemLogConst} forcing $G$ to be \K--Einstein.
\end{proof}

\begin{remark}\rm\label{RmkAlsoOtherMetrics}
	The proofs of Theorems~\ref{mainteor2} and \ref{mainteor3} show that the same conclusions can be achieved taking any metric $g$ such that $D_p^g\in \log\tilde K_p$.
	Examples of such metrics are given by the  \K\ metrics induced by those in Examples~\ref{projindnoncomp}-\ref{exnotwbcomp2} above.
\end{remark}

\bibliographystyle{amsplain}
	
\bibliography{/Users/giovanniplacini/Library/CloudStorage/Dropbox/Projects/RTT_PROGETTI/Bibliography/biblio.bib}

\end{document}